\documentclass[11pt]{amsart}
\usepackage{mathptmx}

\newcommand{\R}{\mathbb{R}}
\newcommand{\Q}{\mathbb{Q}}
\newcommand{\Z}{\mathbb{Z}}

\renewcommand{\a}{\alpha}

\newcommand{\e}{\epsilon}

\newcommand{\exact}[1]{\stackrel{#1}{\longrightarrow}}
\newcommand{\indlim}{\varinjlim}
\newcommand{\prolim}{\varprojlim}
\renewcommand{\hat}{\widehat}
\renewcommand{\epsilon}{\varepsilon}

\newtheorem{thm}{Theorem}

\newtheorem{cor}{Corollary}
\newtheorem{lem}{Lemma}
\newtheorem{prop}{Proposition}

\theoremstyle{definition} %% Body style is Roman.
\newtheorem{Def}{Definition}
\newtheorem{rem}{Remark}

\theoremstyle{remark} %% Headline style is not bold.

\newtheorem*{ack}{Acknowledgment}

\title[Zeros of the Alexander polynomial of knot]{Zeros of the Alexander polynomial of knot}
\author{Akio Noguchi}
\address{
Department of Mathematics,
Tokyo Institute of Technology,
Oh-okayama, Meguro-ku, Tokyo 152-8551, Japan
}
\email{akio@math.titech.ac.jp }
\subjclass[2000]{Primary 57M27; Secondary 11S05, 37B40}
\thanks{The author was supported in part by JSPS fellowship for young scientists.}

\begin{document}

\begin{abstract}
The leading coefficient of the Alexander polynomial of a knot is the most informative element in this invariant, and the growth of orders of the first homology of cyclic branched covering spaces is also a familiar subject. Accordingly, there are a lot of investigations into each subject. However, there is no study which deal with the both subjects in a same context. In this paper, we show that the two subjects are closely related in $p$-adic number theory and dynamical systems.
\end{abstract}

\maketitle

\section{Introduction}
The leading coefficient of the Alexander polynomial $\Delta_K(t)$ of a knot $K$ is a well-known invariant for detecting fibered knots. The Alexander polynomial of a fibered knot is always monic~\cite{MR0176462}. The converse is not always true, but it holds for many knots, for example, alternating knots~\cite{MR0275414}. Moreover the monic condition characterizes fibered knots in a sense of realization~\cite{MR0199858,MR554949}.

The leading coefficient of the Alexander polynomial of a knot is also related to the commutator subgroup $G_K'$ of the knot group $G_K=\pi_1(S^3 \setminus K)$. The abelianaization of $G_K'$ is  finitely generated if and only if the leading coefficient is $\pm1$~\cite{MR27:4226,MR0116047}.

The $r$-fold cyclic covering $X_r(K)$ branched over a knot $K$ is a fundamental object in the knot theory because topological invariants of it are also invariants of the knot. In \cite{MR45:4394}, Gordon studied the growth of the order of $H_1(X_r(K);\Z)$ with respect to $r$ and asked whether the growth is exponential in case some zeros of $\Delta_K(t)$ are not a root of unity. More than 15 years later, this question was affirmatively answered by Riley~\cite{MR92g:57017} and Gonz{\'a}lez-Acu{\~n}a and Short~\cite{MR93g:57004} independently. 

As Gordon commented, the difficulty in computing the growth of orders arises from the case in which all zeros belong to the unit circle but some are not a root of unity (e.g. $\Delta_{5_2}(t)=2t^2-3t+2$). In this case, the standard norm is useless. Reliy~\cite{MR92g:57017} managed the difficulty by the $p$-adic analysis. On the other hand, Gonz{\'a}lez-Acu{\~n}a and Short~\cite{MR93g:57004} managed it by showing that the growth is equal to the Mahler measure of the Alexander polynomial. Moreover, in their introduction, Gonz{\'a}lez-Acu{\~n}a and Short remarked possibility of interpretation of the growth as entropy of a dynamical system.

In this paper, we clear the dynamical system which was predicted by Gonz{\'a}lez-Acu{\~n}a and Short and also explain that the $p$-adic approach of Reliy is useful in essence. 

Solenoidal dynamical systems are nurtural generalizations of toral dynamical systems. The entropy of solenoidal dynamical systems was first calculated by Juzvinski{\u\i}~\cite{MR0214726} and latterly re-calculated by Lind and Ward~\cite{MR90a:28031} with $p$-adic number theory. Furthermore Einseidler and Ward~\cite{MR1777892} also investigate the relation between the Mahler measure and the entropy of a solenoidal dynamical system associated with a Fitting ideal of a module.

Applying the above results of solenoidal dynamical systems to knot theory, we refine some topics of the Alexander polynomial of a knot. Here is a instance of results, in which the apparently different topics are really close. The reason of this similarity are explained in this paper.

Let $\Delta_K(t)=\sum_{i=0}^n a_i t^i \;(a_0a_n\ne0)$ be the Alexander polynomial of a knot $K$ and $\a_{i}$ the zeros (counted with multiplicity) of $\Delta_K(t)$. Then,
\begin{enumerate}
\item the leading coefficient of $\Delta_K(t)$ is
\begin{align*}
\log|a_n|=\sum_{p<\infty} \sum_{|\a_{i}|_p > 1} \log|\a_{i}|_p 
\end{align*}
(Corollary~\ref{a_n}), and
\item the growth of order of the first homology of the $r$-fold cyclic covering branched over $K$ is
\begin{align*}
\lim_{\stackrel{r\to\infty}{|H_1(\cdot)|\ne0}} \frac{\log|H_1(X_r(K);\Z)|}{r}=\sum_{p\le\infty} \sum_{|\a_{i}|_p > 1} \log|\a_{i}|_p 
\end{align*}(Corollary~\ref{growth}).
\end{enumerate}
Here, $|\cdot|_p$ are the $p$-adic norms and $|\cdot|_\infty$ is the standard norm. (we assume that the embedding $\overline{\Q} \to \overline{\Q_p}$ are fixed.)

In our study, we also establish the followings.

\begin{itemize}
\item The distribution of the zeros measures a ``distance" of the Alexander module from being finitely generated as $\Z$-module. (Section 3.2)
\item In~\cite{MR45:4394}, the primary interest in investigating the growth (2) might be to study the periodicity of $H_1(X_r(K);\Z)$. However, the growth also measures complexity of the Alexander module naturally in a sense of solenoidal dynamical systems. This means, this growth can be interpreted as volume growth in an adele ring. (Section 3.1)
\end{itemize}

\smallskip

Here, we have a few comments on this study. It might make our study a little more attractive.

~

The Alexander polynomial of a knot is defined as a greatest common divisor of the initial Fitting ideal (elementary ideal) of the Alexander module $H_1(X_\infty(K);\Z)$ as $\Z[t^\pm]$-module. Here, the indeterminate $t$ is identified with the meridian action on $H_1(X_\infty(K);\Z)$ and $X_\infty(K)$ is the infinite cyclic cover of $X(K)=S^3 \setminus K$. Then, by tensoring with the rational numbers $\Q$, the Alexander polynomial is also generator of the Fitting ideal of the module $H_1(X_\infty(K);\Q)$ as $\Q[t^\pm]$-module, and hence it is the characteristic polynomial of the meridian action on $H_1(X_\infty(K);\Q)$, up to units (see Theorem 6.17 in \cite{MR98f:57015}).

While the rational homology $H_1(X_\infty(K);\Q)$ gives a nice explanation of the Alexander polynomial, the leading coefficient $a_n$ is lost in $H_1(X_\infty(K);\Q)$ because $a_n$ is a unit in $\Q[t^\pm]$. On the other hand, the entropy has an advantage over the Fitting ideal because we can replace $H_1(X_\infty(K);\Z)$ with $H_1(X_\infty(K);\Q)$ with preserving the entropy (cf. Step 1 in the proof of Proposition~\ref{LW}). %This is why the zeros of the Alexander polynomial keep information about $H_1(X_\infty(K);\Z)$.

~

In~\cite{MR92g:57017}, Riley also proved a result on the $p$-part of $|H_1(X_r(K);\Z)|$. He obtained the upper bounds for the $p$-parts: $|H_1(X_r(K);\Z)|^{(p)}< A H^p{^E} r^n$, where $A, H, E, n$ are constants depending on a knot and this upper bounds is best possible (up to constant multiplier). That is, the every $p$-parts growth at most polynomially with respect to $r$ nevertheless the entire part $|H_1(X_r(K);\Z)|$ can growth exponentially. 
(Silver and Williams~\cite{MR1955605} also studied this topic. Unfortunately, this part of their results had been already established by Riley, however they discussed this topic under mild hypotheses.)

~

In this paper, the $p$-adic coefficient homology group $H_1(X_\infty(K);\Q_p)$ plays an important role in studying the distribution of zeros (Theorem 1). This approach was motivated by the previous paper~\cite{math.GT/0505531}, in which we discussed that an analogy between the reciprocity of $\Delta_K$ and the functional equation in the Weil conjecture. Then, the $p$-adic (co)homology theory is also investigated to approach the Weil conjecture.

\begin{ack}
I would like to thank Professor Kazuo Masuda for his helpful advice and valuable discussion, and also thank Professor Sadayoshi Kojima, Professor Hitoshi Murakami, Professor Masanori Morishita and Professor Gregor Masbaum for their kind reading my rough draft and giving helpful comments. And, I am grateful to Professor Kunio Murasugi for his encouragement in this study and kind hospitality during my visiting University of Toronto. I also appreciate a seminor with Professor Shoichi Nakajima, Professor Shin Nakano and Professor Mikami Hirasawa, which was helpful to meke this paper readable.
\end{ack}

\section{Solenoidal entropy and Alexander polynomial}
A solenoid $\Sigma^n$ is, by definition, a compact connected finite-dimensional abelian group, which was arose from a generalization of the torus $\mathbb{T}^n$. The following theorem, which was given by Lind and Ward~\cite{MR90a:28031}, plays a key role for our results.

\begin{prop}[Lind and Ward~\cite{MR90a:28031}]\label{LW}
Let $T$ be an automorphism of $n$-dimensional solenoid $\Sigma^n$. Then,
\begin{enumerate}
\item the entropy of $T$ is the sum of the entropies of the automorphisms of $\Q_p^n$ induced by $T$
\begin{align*}
h(T;\Sigma^n)=\sum_{p\le\infty} h(T;\Q_p^n),
\end{align*}
and
\item the $p$-adic entropy is computed by the eigenvalues $\lambda_{1},\dots,\lambda_{n}$ of the induced automorphism in $GL(d,\Q_n)$ as follows:
\begin{align*}
h(T;\Q_p^n)=\sum_{|\lambda_{i}|_p > 1} \log|\lambda_{i}|_p .
\end{align*}
\end{enumerate}
\end{prop}

We shall confuse the topological entropy and the measure theoretic entropy by the formula of Bowen~\cite{MR43:469}:
\begin{align*}
h_d(T)=\lim_{\e \to 0} \limsup_{n \to \infty} [ -\frac{1}{n} \log \mu(\cap_{k=0}^{n-1} T^{-k}B(e,\e)) ],
\end{align*}
where $B(e,\e)$ is a open $\e$-ball of the identity element with respect to an invariant metric $d$, $\mu$ is a Haar measure and $T$ a surjective endomorphism on a locally compact abelian group. In other word, we regard the entropy as the topological entropy and also the measure theoretic entropy with respect to a Haar measure.

In \cite{MR90a:28031}, Lind and Ward computed the entropy of solenoidal automorphisms with intrinsic arguments. The proof is really helpful for our applications later. Although we review an outline of the proof, referring to the original paper is strongly recommend.

\begin{proof}[Outline of the proof]
~
\begin{enumerate}
\item [Step 1:] Because the dual group of $\Sigma^d$ can be embed into $\Q^d$ $\indlim\Gamma_n \cong \Q^d$ for $\Gamma_n=\frac{1}{n!}\hat{\Sigma}^d$. That is, $\hat{\Q}^d\cong\prolim \hat{\Gamma}_n$. Hence there are $K_n$ such that $\hat{\Gamma}_n\cong\hat{\Q}^d/K_n$ and consequently $h(T;\hat{\Q}^d)=h(T;\hat{\Gamma}_n) + h(T;K_n)$. Because $h(T;\hat{\Gamma}_n)=h(T;\Sigma^d)$ for any $n$ and $h(T;K_n) \to 0$ as $n\to\infty$, we obtain that
\begin{align*}
h(T;\Sigma^d)=h(T;\widehat{\Q}^d).
\end{align*}
\item [Step 2:] The entropy on the full solenoid $\hat{\Q}^d$ can be lifted to the entropy on the adele ring $\mathbb{A}_\Q^d$ because $\mathbb{A}_\Q/\Q \cong \hat{\Q}$. So, we have
\begin{align*}
h(T;\hat{\Q}^d)=h(T;\mathbb{A}_{\Q}^d).
\end{align*}
\item [Step 3:] Since the adele ring is a restricted direct product space, the entropy on it can be decomposed into the entropies of each direction, i.e.
\begin{align*}
h(T;\mathbb{A}_{\Q}^d)=\sum_{p\le\infty}h(T;\Q_p^d).
\end{align*}
\item [Step 4:] Finally, the $p$-adic norm is the module on $\Q_p$ (i.e. $\mu(a\cdot)=|a|_p\mu(\cdot)$ for any $a\in\Q_p$, where $\mu$ is a Haar measure), we have
\begin{align*}
h(T;\Q_p^d)=\sum_{|\lambda_{i}|_p > 1} \log|\lambda_{i}|_p.
\end{align*}
\end{enumerate}
\end{proof}

To connect the Alexander polynomial with a solenoidal automorphism, we need the following lemma.
\begin{lem}\label{sol}
For any knot, the dual group of the first homology group of an infinite cyclic cover $H_1(X_\infty(K);\Z)$ is a $n$-dimensional solenoid. Here, $n$ is the degree of the Alexander polynomial of the knot.
\end{lem}

\begin{proof}
It is sufficient to prove that $H_1(X_\infty(K);\Z)$ is a discrete torsion-free abelian group which has finite-rank $n$. Rapaport~\cite{MR0116047} and Crowell~\cite{MR27:4226} proved that $H_1(X_\infty(K);\Z)$ is torsion-free and has finite rank $n$. (Here, the rank of $A$ means the cardinality of any maximal set of $\Z$-linearly independent elements of $A$.)
\end{proof}

Because the Alexander polynomial is equal to the characteristic polynomial of the  meridian action on $H_1(X_\infty(K);\Q)$, up to multiplication by a unit, the following theorem follows immediately from Proposition~\ref{LW} and Lemma~\ref{sol}.

\begin{thm}\label{thm}
Let $\a_{i}$ be the zeros (counted with multiplicity) of the Alexander polynomial of a knot. Then,
\begin{enumerate}
\item the entropy of the meridian action on the $p$-adic Alexander module;
$t_p:H_1(X_\infty(K);\Q_p) \to H_1(X_\infty(K);\Q_p)$ is
\begin{align*}
h(t_p)=\sum_{|\a_{i}|_p > 1} \log|\a_{i}|_p ,
\end{align*}
where $|\cdot|_p$ is the $p$-adic norms, and
\item the entropy of dual action of meridian $\hat{t}:\widehat{H_1(X_\infty(K);\Z)} \to \widehat{H_1(X_\infty(K);\Z)}$ is $h(\hat{t})=\sum_{p\le\infty} h(t_p)$, that is
\begin{align*}
h(\hat{t})=\sum_{p\le\infty} \sum_{|\a_{i}|_p > 1} \log|\a_{i}|_p.
\end{align*}
\end{enumerate}
Here $\Q_\infty = \R$ by convention.
\end{thm}

\section{Applications}

\subsection{Growth of order of homology of branched cyclic covering space}

In this section, we study the relation between the $p$-adic norm of the zeros of $\Delta_K$ and the growth of orders of the first homology groups of the $r$-fold cyclic covering of $S^3$ branched over $K$. Roughly speaking, the $r$-fold cyclic covering branched over $K$ is a compact space which associated with a homomorphism $G=\pi_1(S^3 \setminus K) \to \Z/r\Z$ (for precise definition, see \cite{MR87b:57004,MR98f:57015}).

The order of the first homology group of this space can be computed by the following formula.

\begin{prop}[Fox~\cite{MR0095876}, Weber~\cite{MR570312}]\label{Fox}
Let $X_r(K)$ be the $r$-fold cyclic covering of $S^3$ branched over $K$. Then, the order of the first homology group of $X_r(K)$ is given by
$$
  \left|H_1(X_r(K);\Z)\right|
  =
  \left|\prod_{d=1}^{r-1}\Delta_K(\exp({2d\pi\sqrt{-1}}/{r}))\right|.
$$
By convention, $|H_1(X_r(K);\Z)|=0$ means that $H_1(X_r(K);\Z)$ is an infinite group.
\end{prop}

\begin{Def}[logarithmic Mahler measure~\cite{MR0124467}]
For non-zero Laurent polynomial $f(x)$ with integral coefficients, the logarithmic Mahler measure of $f$ is defined by
\begin{align*}
m(f)=\int^1_0 \log|f(\exp(2\pi t \sqrt{-1}))|dt.
\end{align*}
\end{Def}

The growth of orders $\left|H_1(X_r(K);\Z)\right|$ is expressed by the logarithmic Mahler measure of the Alexander polynomial, which was proved by Gonz\'alez-Acu\~na and Short~\cite{MR93g:57004}.

\begin{prop}[Gonz\'alez-Acu\~na and Short~\cite{MR93g:57004}]\label{GAS}
\begin{align*}
\lim_{\stackrel{r\to\infty}{|H_1(\cdot)|\ne0}} \frac{\log|H_1(X_r(K);\Z)|}{r} = m(\Delta_K).
\end{align*}
\end{prop}

\begin{rem}
In~\cite{MR2003h:57011}, Silver and Williams generalized the result of Gonz\'alez-Acu\~na and Short~\cite{MR93g:57004} to links and also pointed out that Proposition~\ref{GAS} holds though the order $|H_1(X_r(K);\Z)|$ is replaced with the order of the torsion subgroup $|TH_1(X_r(K);\Z)|$.
%with dynamical systems approach which was predicted by Gonz\'alez-Acu\~na and Short, and our study is motivated by their work.
\end{rem}

By the way, the Mahler measure is deeply related to the entropy of an algebraic dynamical system, which is found in~\cite{MR92j:22013} for example. In this paper, we use a more suitable result which was proved by Einseidler and Ward~\cite{MR1777892}.

\begin{prop}[Einseidler and Ward~\cite{MR1777892}]\label{EW}
Let $ 0 \exact{} F_n \exact{} \dots \exact{} F_1 \exact{\phi_1} F_0 \exact{} M \exact{} 0$ be a finite free resolution of the $\Z[t^\pm]$-module $M$ and $J(\phi_1)$ the initial Fitting ideal. Let $\a_t$ a natural automorphism which is induced by the shift of the indeterminate $t$. Then the entropy of $\a_t$ is
$$h(\hat{\a_t})=m(\gcd(J(\phi_1))).$$
\end{prop}

By combining Proposition~\ref{Fox} and~\ref{EW}, we can see that the growth of orders gives another method to compute the entropy in Theorem~\ref{thm}. Therefore, we can obtain the following corollary.

\begin{cor}\label{growth}
\begin{align*}
\lim_{\stackrel{r\to\infty}{|H_1(\cdot)|\ne0}} \frac{\log|H_1(X_r(K);\Z)|}{r} 
  = \sum_{p\le\infty}\sum_{|\a_{i}|_p > 1} \log|\a_{i}|_p, \label{infsum} 
\end{align*}
where $\a_{i}$ are the zeros of the Alexander polynomial $\Delta_K(t)$.
\end{cor}

Proposition \ref{GAS} and Corollary \ref{growth} are similar. However, the every terms of the right hand side of Corollary~\ref{growth} make sense in a dynamical system. Actually, Corollary~\ref{growth} can be rewritten into the following form:
\begin{align*}
\lim_{\stackrel{r\to\infty}{|H_1(\cdot)|\ne0}} \frac{\log|H_1(X_r(K);\Z)|}{r} 
  = \sum_{p\le\infty} h(t_p),
\end{align*}
where $h(t_p)$ is the entropy of the meridian action $t_p:H_1(X_\infty(K);\Q_p) \to H_1(X_\infty(K);\Q_p)$. In other words, we resolved the growth into fine factors of the $p$-adic entropy.

Needless to say, we also obtain the following property as the spacial case of Corollary~\ref{growth}.

\begin{cor}[Riley~\cite{MR92g:57017}, Gonz\'alez-Acu\~na and Short~\cite{MR93g:57004}]\label{units}
Let $X_r(K)$ be the $r$-fold cyclic covering branched over $K$.
Then, if the Alexander polynomial $\Delta_K(t)$ have zeros which are not roots of unity,
the finite values of the order of the first homology group $|H_1(X_r(K);\Z)|$ grows exponentially with respect to $r$.
\end{cor}

\begin{proof}
(Indirect proof) Because all $\a_{i}$ belong to the valuation ring $\mathcal{O}_p=\{x\in \overline{\Q_p} \mid |x|_p \le 1 \}$, $f(t)=\prod_i(t-\a_i)$, which is $\Delta(t)$ up to constant multiples, belongs to $\Z_p[t]\cap\Q[t]$, where $\Z_p=\{x\in\Q_p \mid |x|_p\le1 \}$. This holds for any prime $p$. Hence $f(t)\in\Z[t]$. (Another way to see this is to prove Corollary~\ref{a_n} at the first. But in this proof, the condition $\Delta(1)=\pm1$ is not necessary.) Consequently, the zeros of the Alexander polynomial must be roots of unity from $|\a_i|\le1$ and the Kronecker's theorem~\cite{Kronecker}.
\end{proof}

\subsection{Leading coefficient of Alexander polynomial}

In this section, we apply Theorem 1 to a criterion for being finitely generated as $\Z$-module. %Formerly, we have utilized the leading coefficient of the Alexander polynomial for the criterion. But it only determines whether the Alexander module is finitely generated or not. Now, we reveal that the leading coefficient of the Alexander polynomial is not just the criterion for being finitely generated as $\Z$-module.
By the corollary below, we can regard the entropies $h(t_p)$ for all primes $p<\infty$ as obstructions for being finitely generated. 

\begin{cor}\label{criterion}
Let $h(t_p)$ be the entropy of the meridian action on the $p$-adic Alexander module $H_1(X_\infty(K);\Q_p)$. 
If the Alexander module $H_1(X_\infty(K);\Z)$ is finitely generated as $\Z$-module, then all the entropies $h(t_p)$ are equal to zero for all finite primes $p<\infty$.
\end{cor}

\begin{proof}
The Alexander module $H_1(X_\infty(K);\Z)$ is finitely generated if and only if $\hat{H_1(X_\infty(K);\Z)}$ is isomorphic to the $n$ dimensional torus. Then, it is covered by $H_1(X_\infty(K);\R)$. By the well-known result for toral automorphisms, the entropy of the meridian action on $\hat{H_1(X_\infty(K);\Z)}$ is
\begin{align*}
h(t)=\sum_{|\a_i| > 1} \log|\a_i|,
\end{align*}
where $\a_i$ are eigenvalues of the meridian action on $H_1(X_\infty(K);\R)$. 

These entropy must be equal to the entropy in Theorem 1-(2). Therefore, the entropies of the meridian action on $H_1(X_\infty(K);\Q_p)$ are zero for any $p<\infty$, that is
$$h(t_p)=\sum_{|\a_{i}|_p > 1} \log|\a_{i}|_p=0 \quad \text{for} \quad p<\infty.$$
\end{proof}

This obstructions give a new viewpoint of the leading coefficient of the Alexander polynomial. The following Corollary means that the entropies $h(t_p)$ are fine factors of the leading coefficient of the Alexander polynomial.

\begin{cor}\label{a_n}
Let $\a_{i}$ be the zeros of $\Delta_K(t)=\sum_{i=0}^n a_i t^i$. Then the leading coefficient of $\Delta_K(t)$ is the sum of the entropies of the meridian action on the $p$-adic Alexander module $H_1(X_\infty(K);\Q_p)$ for the finite primes $p<\infty$. In other words, the following equation holds:
\begin{align*}
\log |a_n| &= \sum_{p<\infty}\sum_{|\a_{i}|_p > 1} \log|\a_{i}|_p.
\end{align*}
\end{cor}

\begin{proof}
Let $f(t)=\Delta_K(t)/a_n=\prod(t-\a_{i})$ and $s$ the least common multiple of the denominators of the coefficients of $f(t)$. Then, 
\begin{align*}
\sum_{p<\infty} h(t_p) = \sum_{p<\infty}\sum_{|\a_{i}|_p > 1} \log|\a_{i}|_p = \log s.
\end{align*}
Because $\Delta_K(1)=\pm 1$, the coefficients are relatively prime $(a_n,\dots,a_1)=1$. Hence $s=|a_n|$.
(The above argument is essentially found in the proofs of Theorem 3 in \cite{MR90a:28031} and Theorem 2 in \cite{MR92g:57017}.)
\end{proof}

Since $\Delta_K(1)=\pm1$ for any knot, the Alexander polynomial can be completely determined (up to $\pm1$) by the zeros. Hence the leading coefficient is also determined. On the other hand, Corollary~\ref{a_n} shows that the leading coefficient is resolved into the $p$-adic entropy and consequently can be recovered from the zeros.

\subsection{Final remarks}

\subsubsection{Determining knots by cyclic branched covers}
In~\cite{MR903865}, Kojima showed that prime knots are determined by their cyclic branched covers. So, there might be a method to compute the Alexander module by the data of cyclic branched covers.

The growth of orders $|H_1(X_r(K);\Z)|$ does not determine completely the Alexander module $H_1(X_\infty(K);\Z)$. However the growth measures complexity of the Alexander module from dynamical viewpoint, and also it is a similar invariant with the leading coefficient of the Alexander polynomial.

In addition, it is still remain open whether infinitely many branched covers are necessary for determining the knot. %We expect the stronger method to reconstruct the Alexander module or other invariants by the data of cyclic branched covers.

\subsubsection{Volume conjecture}
The volume (or Kashaev) conjecture~\cite{MR1434238,MR1828373} expects that the asymptotic behavior of the Kashaev invariant (= the specialization of the colored Jones polynomial) implies the hyperbolic volume of the complement of the knot (hyperbolic case).

In general, the topological entropy picks up natural measures by the variational principle. In other words, the topological entropy is related with natural measure theory behind it.

In this paper, the asymptotic behavior of the special values of the Alexander polynomial is related with the entropy of an action on the adele ring of $\Q$. Consequently, the asymptotic behavior is related with the volume growth with respect to the Haar measure on the adele ring.

% Then, the topological entropy is equal to the measure theoretic entropy with respect to a measure of the maximal entropy. In our case, the Haar measure on the Alexander module (precisely, the Plancherel measure on its dual group) is a measure of the maximal entropy, and this measure is lifted to adele rings with preserving the entropy.

%Consequently the asymptotic behavior (growth of orders) can be translated into expansions of volumes with respect to the Haar measure on the adele ring because the growth is interpreted as the entropy of the meridian action.

Our point is that, when an asymptotic behavior is interpreted as an entropy, it can be related with a natural measure theory. Does this strategy work out for the volume conjecture?


\begin{thebibliography}{10}

\bibitem{MR43:469}
R.~Bowen.
\newblock Entropy for group endomorphisms and homogeneous spaces.
\newblock {\em Trans. Amer. Math. Soc.}, 153:401--414, 1971.

\bibitem{MR0199858}
G.~Burde.
\newblock Alexanderpolynome {N}euwirthscher {K}noten.
\newblock {\em Topology}, 5:321--330, 1966.

\bibitem{MR87b:57004}
G.~Burde and H.~Zieschang.
\newblock {\em Knots}, volume~5 of {\em de Gruyter Studies in Mathematics}.
\newblock Walter de Gruyter \& Co., Berlin, 1985.

\bibitem{MR27:4226}
R.~H. Crowell.
\newblock The group {$G'/G''$} of a knot group {$G$}.
\newblock {\em Duke Math. J.}, 30:349--354, 1963.

\bibitem{MR1777892}
M.~Einsiedler and T.~Ward.
\newblock Fitting ideals for finitely presented algebraic dynamical systems.
\newblock {\em Aequationes Math.}, 60(1-2):57--71, 2000.

\bibitem{MR0095876}
R.~H. Fox.
\newblock Free differential calculus. {III}. {S}ubgroups.
\newblock {\em Ann. of Math. (2)}, 64:407--419, 1956.

\bibitem{MR93g:57004}
F.~Gonz{\'a}lez-Acu{\~n}a and H.~Short.
\newblock Cyclic branched coverings of knots and homology spheres.
\newblock {\em Rev. Mat. Univ. Complut. Madrid}, 4(1):97--120, 1991.

\bibitem{MR45:4394}
C.~M. Gordon.
\newblock Knots whose branched cyclic coverings have periodic homology.
\newblock {\em Trans. Amer. Math. Soc.}, 168:357--370, 1972.

\bibitem{MR0214726}
S.~A. Juzvinski{\u\i}.
\newblock Computing the entropy of a group of endomorphisms.
\newblock {\em Sibirsk. Mat. \u Z.}, 8:230--239, 1967.
\newblock English transl. in Siberian Math. Journal 8 (1968), 172-178.

\bibitem{MR1434238}
R.~M. Kashaev.
\newblock The hyperbolic volume of knots from the quantum dilogarithm.
\newblock {\em Lett. Math. Phys.}, 39(3):269--275, 1997.

\bibitem{MR903865}
S.~Kojima.
\newblock Determining knots by branched covers.
\newblock In {\em Low-dimensional topology and Kleinian groups
  (Coventry/Durham, 1984)}, volume 112 of {\em London Math. Soc. Lecture Note
  Ser.}, pages 193--207. Cambridge Univ. Press, Cambridge, 1986.

\bibitem{Kronecker}
L.~Kronecker.
\newblock {Z}wei {S}\"atze ueber {G}leichungen mit ganzzahligen
  {C}oeffichienten.
\newblock {\em J. Reine Angew. Math.}, 53:173--175, 1857.

\bibitem{MR98f:57015}
W.~B.~R. Lickorish.
\newblock {\em An introduction to knot theory}, volume 175 of {\em Graduate
  Texts in Mathematics}.
\newblock Springer-Verlag, New York, 1997.

\bibitem{MR92j:22013}
D.~Lind, K.~Schmidt, and T.~Ward.
\newblock Mahler measure and entropy for commuting automorphisms of compact
  groups.
\newblock {\em Invent. Math.}, 101(3):593--629, 1990.

\bibitem{MR90a:28031}
D.~A. Lind and T.~Ward.
\newblock Automorphisms of solenoids and {$p$}-adic entropy.
\newblock {\em Ergodic Theory Dynam. Systems}, 8(3):411--419, 1988.

\bibitem{MR0124467}
K.~Mahler.
\newblock An application of {J}ensen's formula to polynomials.
\newblock {\em Mathematika}, 7:98--100, 1960.

\bibitem{MR1828373}
H.~Murakami and J.~Murakami.
\newblock The colored {J}ones polynomials and the simplicial volume of a knot.
\newblock {\em Acta Math.}, 186(1):85--104, 2001.

\bibitem{MR0275414}
K.~Murasugi.
\newblock The commutator subgroups of the alternating knot groups.
\newblock {\em Proc. Amer. Math. Soc.}, 28:237--241, 1971.

\bibitem{MR0176462}
L.~P. Neuwirth.
\newblock {\em Knot groups}.
\newblock Annals of Mathematics Studies, No. 56. Princeton University Press,
  Princeton, N.J., 1965.

\bibitem{math.GT/0505531}
A.~Noguchi.
\newblock A functional equation for the lefschetz zeta functions of infinite
  cyclic coverings with an application to knot theory.
\newblock preprint, arXiv:math.GT/0505531.

\bibitem{MR554949}
C.~V. Quach.
\newblock Polyn\^ome d'{A}lexander des noeuds fibr\'es.
\newblock {\em C. R. Acad. Sci. Paris S\'er. A-B}, 289(6):A375--A377, 1979.

\bibitem{MR0116047}
E.~S. Rapaport.
\newblock On the commutator subgroup of a knot group.
\newblock {\em Ann. of Math. (2)}, 71:157--162, 1960.

\bibitem{MR92g:57017}
R.~Riley.
\newblock Growth of order of homology of cyclic branched covers of knots.
\newblock {\em Bull. London Math. Soc.}, 22(3):287--297, 1990.

\bibitem{MR2003h:57011}
D.~S. Silver and S.~G. Williams.
\newblock Mahler measure, links and homology growth.
\newblock {\em Topology}, 41(5):979--991, 2002.

\bibitem{MR1955605}
D.~S. Silver and S.~G. Williams.
\newblock Torsion numbers of augmented groups with applications to knots and
  links.
\newblock {\em Enseign. Math. (2)}, 48(3-4):317--343, 2002.

\bibitem{MR570312}
C.~Weber.
\newblock Sur une formule de {R}. {H}. {F}ox concernant l'homologie des
  rev\^etements cycliques.
\newblock {\em Enseign. Math. (2)}, 25(3-4):261--272 (1980), 1979.

\end{thebibliography}
\end{document}